\documentclass[a4paper]{amsart}
\usepackage{amsmath}
\usepackage{amssymb}
\usepackage{amsthm}
\usepackage{amscd}

\usepackage{ulem}
\usepackage{comment}

\usepackage{xcolor}
\usepackage{hyperref}

\newenvironment{solution}{\begin{proof}[Proof.]\parskip=0pt}{\end{proof}}

\newtheorem*{theorem}{Theorem}
\newtheorem*{proposition}{Proposition}
\newtheorem{thm}{Theorem}[section]
\newtheorem{lemma}[thm]{Lemma}
\newtheorem{prop}[thm]{Proposition}

\theoremstyle{definition}

\newtheorem{expl}[thm]{Example}

\theoremstyle{remark}
\newtheorem*{rem}{Remark}

\numberwithin{equation}{section}

\def\k{\mathbf{k}}
\def\K{\mathcal{K}}
\def\H{\mathcal{H}}
\def\Z{\mathcal{Z}}

\def\D{\Delta}

\def\x{\times}
\def\ox{\otimes}

\def\sm{\setminus}
\def\ss{\subset}

\def\wt{\widetilde}

\def\mb{\mathbb}

\def\oname{\operatorname}

\newcommand{\6}{\partial}

\hypersetup{
    colorlinks=true,
    linkcolor=blue,
    urlcolor=blue,
    citecolor=red,
}
\title[Moment-angle manifolds, chordality and connected sums]{Moment-angle manifolds corresponding to three-dimensional simplicial spheres, chordality and connected sums of products of spheres}
\author{Victoria Kovyrshina}
\address{Department of Mathematics and Mechanics, Moscow
State University, Russia;\newline
National Research University Higher School of Economics, Moscow, Russia}
\email{potchtovy\_jashik@mail.ru}

\author{Taras Panov}
\address{
Steklov Mathematical Institute of Russian Academy of Sciences, Moscow
\newline
Department of Mathematics and Mechanics, Moscow
State University, Russia;\newline
National Research University Higher School of Economics, Moscow, Russia}
\email{tpanov@mi-ras.ru}

\thanks{This work was supported by the Russian Science Foundation under grant no.~24-71-00059, https://rscf.ru/project/24-71-00059/. Victoria Kovyrshina is a recipient of a scholarship from the Theoretical Physics and Mathematics Advancement Foundation ``BASIS''}

\begin{document}
\begin{abstract}
We prove that the moment-angle complex $\Z_\K$ corresponding to a $3$-dimensional simplicial sphere $\K$ has the cohomology ring isomorphic to the cohomology ring of a connected sum of products of spheres if and only if either (a) $\K$ is the boundary of a $4$-dimensional cross-polytope, or (b) the one-skeleton of $\K$ is a chordal graph, or (c) there are only two missing edges in $\K$ and they form a chordless $4$-cycle. For simplicial spheres $\K$ of arbitrary dimension, we obtain a sufficient condition for the ring isomorphism $H^*(\mathcal Z_\K)\cong H^*(M)$ where $M$ is a connected sum of products of spheres.
\end{abstract}

\maketitle

\section{Introduction}

The moment-angle complex is a topological space (a CW complex) with a torus action that features in toric topology and homotopy theory of polyhedral products~\cite{BP}. The topology of a moment-angle complex $\Z_\K$ is determined by the combinatorics of the corresponding simplicial complex~$\K$.
If $\K$ is the nerve complex of a simple polytope $P$, then the corresponding moment-angle complex, which is denoted by $\Z_P$, is a smooth manifold. 

There are several different geometric constructions of moment-angle manifolds enriching their topology with remarkable and peculiar geometric structures. One of them arises in holomorphic dynamics, where the moment-angle manifold $\Z_P$ appears as the leaf space of a holomorphic foliation on an open subset of a complex space, and is diffeomorphic to a nondegenerate intersection of Hermitian quadrics~\cite{bo-me06}, \cite[Chapter~6]{BP}. All early examples of moment-angle manifolds appearing in this context where diffeomorphic to connected sums of products of spheres. This is the case, for example, when $P$ is two-dimensional (a polygon). From the description of the cohomology ring of $\Z_P$ it became clear that the topology of moment-angle manifolds in general is much more complicated than that of a connected sum of sphere products; for instance, $H^*(\Z_P)$ can have arbitrary additive torsion or nontrivial higher Massey products~\cite[Chapter~4]{BP}. 

Nevertheless, the question remained of identifying the class of simple polytopes $P$ (or more generally, simplicial spheres~$\K$) for which the moment-angle manifold $\Z_P$ is homeomorphic to a connected sum of products of spheres. This question is also interesting from the combinatorial and homotopy-theoretic points of view, as it is related to the conditions for the minimal non-Golodness of $\mathcal K$ and the chordality of its one-skeleton. For three-dimensional polytopes $P$ (or two-dimensional spheres $\K$), it was proved in~\cite[Proposition~11.6]{bo-me06} that $\Z_P$ is diffeomorphic to a connected sum of products of spheres if and only if $P$ is obtained from a 3-simplex by consecutively cutting off some $l$ vertices. This characterisation can be extended by adding two more equivalent conditions, the chordality and the minimal non-Golodness (see Proposition~\ref{mainthm3}):

\begin{proposition}
Let $\K$ be a two-dimensional simplicial sphere and let $P$ be the a three-dimensional simple polytope such that $\K=\K_P$. Suppose that $P$ is not a cube. The following conditions are equivalent:
\begin{itemize}
    \item[(a)] $P$ is obtained from a simplex $\D^3$ by iterating the vertex cut operation, i.\,e. $P^*$ is a stacked polytope;
    \item[(b)] $\Z_P$ is diffeomorphic to a connected sum of products of spheres;
    \item[(c)] $H^*(\Z_P)$ is isomorphic to the cohomology ring of a connected sum of products of spheres;
    \item[(d)] the one-dimensional skeleton of the nerve complex $\K_P$ is a chordal graph;
    \item[(e)] $\K_P$ is minimally non-Golod, unless $P=\Delta^3$.
\end{itemize}
\end{proposition}

For three-dimensional simplicial spheres $\K$ (including the nerve complexes of four-dimensional simple polytopes) we characterise the moment-angle manifolds $\Z_K$ with the cohomology ring isomorphic to the cohomology ring of a connected sum of products of spheres (see Theorem~\ref{chordal4iff}):

\begin{theorem}
Let $\K$ be a three-dimensional simplicial sphere. There is a ring isomorphism $H^*(\Z_\K) \cong H^*(M_1\#\cdots\# M_k)$ where each $M_i$ is a product of spheres if and only if one of the following conditions is satisfied:
\begin{itemize}
\item[(a)] $\K=S^0*S^0*S^0*S^0$ (the boundary of a $4$-dimensional cross-polytope);
\item[(b)] $\K^1$ is a chordal graph;
\item[(c)] $\K^1$ has exactly two missing edges which form a chordless $4$-cycle.
\end{itemize}
\end{theorem}

We conjecture that under each of the conditions (b) and (c) above the moment-angle manifold $\Z_\K$ is homeomorphic to a connected sum of products of spheres. Under condition (c) we have $H^*(\Z_\K) \cong H^*(M_1\#\cdots\# M_k)$ where one of the summands $M_i$ is a product of \emph{three} spheres. The first example of such $\K$ was constucted in~\cite{FCMW}.

When $\dim P \ge 5$, the chordality of $\K^1_P$ does not imply that $H^*(\Z_P) \cong H^*(M)$ where $M$ is a connected sum of products of spheres, see Example~\ref{5dimex}. A stronger sufficient condition valid for simplicial spheres of arbitrary dimension is given in Theorem~\ref{ssndim}.

\section{Preliminaries}

Let $\K$ be a simplicial complex on the set $[m]=\{1, \ldots, m\}$. We assume that $\K$ contains an empty set $\varnothing$ and all one element subsets $\{i\} \ss [m]$. The dimension of a simplicial complex $\K$ is the maximal cardinality of its simplices minus one.

We denote the full subcomplex of $\K$ on a vertex set $J = \{j_1, \ldots,j_k\} \ss [m]$ by $\K_J$ or by $\K_{\{j_1, \ldots, j_k\}}$.


The \emph{moment-angle complex} $\Z_\K$ corresponding to $\K$ is defined as follows (see \cite[\S4.1]{BP}):
\[
\Z_\K = \underset{I \ss \K}{\bigcup} \Bigl( \underset{i \in I}{\prod} D^2 \x \underset{i \notin I}{\prod} S^1 \Bigr)
\subset\prod_{i=1}^m D^2\,.
\]

\begin{lemma}\label{retract}
If $\K_J$ is a full subcomplex of $\K$, then $\Z_{\K_J}$ is a retract of $\Z_\K$, and $H^*(\Z_{\K_J})$ is a subring of $H^*(\Z_\K)$.
\end{lemma}
\begin{solution}
Let $i\colon \Z_\K \hookrightarrow (D^2)^m$ be the canonical inclusion, and let $q\colon (D^2)^m \rightarrow (D^2)^{|J|}$ be the map that omits the coordinates corresponding to $[m] \sm J$. Then $r = q \circ i\colon \Z_\K \rightarrow \Z_{\K_J}$ is the required retraction, and it induces an injective homomorphism $H^*(\Z_{\K_J})\to H^*(\Z_\K)$ in cohomology.
\end{solution}




\begin{thm}[{\cite[Theorem 4.5.8]{BP}}]\label{Hochster}
    There are isomorphisms of groups 
\[
H^l(\Z_\K) \cong \underset{J \ss [m]}{\bigoplus} \wt{H}^{l - |J| - 1}(\K_J) 
\]
    These isomorphisms combine to form a ring isomorphism $H^*(\Z_\K) \cong \underset{J \ss [m]}{\bigoplus} \wt{H}^*(\K_J)$, where the ring structure on the right hand side is given by the canonical maps
\[
H^{k - |I| - 1}(\K_I)  \ox H^{l - |J| - 1}(\K_J) \longrightarrow H^{k + l - |I| - |J| - 1}(\K_{I \cup J}) \,,
\]
    which are induced by simplicial maps $\K_{I \cup J} \rightarrow \K_I * \K_J$ for $I \cap J = \varnothing$ and zero otherwise.
\end{thm}


We denote 
\[
  \H^{l, J} = \wt{H}^l(\K_J),\quad 
  \H^{*, J} = \wt{H}^*(\K_J)\quad\text{and}\quad
  \H^{l, *} = \underset{J \ss [m]}{\bigoplus} \wt{H}^l(\K_J).
\]  
The ring structure in $H^*(\Z_\K)=\H^{*,*}(\K)$ is given by the maps
\begin{equation}\label{multH}
  \H^{k,I}\otimes\H^{l,J}\longrightarrow \H^{k+l+1,I\sqcup J},\qquad k,l\ge0,\; I\cap J=\varnothing.
\end{equation}

\begin{prop}\label{cohom_l}
If $\K$ is a simplicial complex of dimension $n-1$, then the cohomological product length of $\Z_\K$ is at most~$n$.
\end{prop}
\begin{solution}
Suppose there are elements $c_i \in H^{l_i}(\Z_\K)$, $i=1,\ldots,r$, such that
$c_1 \cdots c_r = c \neq 0$.
This implies, by Theorem~\ref{Hochster}, that there are elements
$\widehat{c} \in \wt{H}^{l}(\K_{J})$ and $\widehat{c}_i \in \wt{H}^{l_i - |J_i| - 1}(\K_{J_i})$ such that $\widehat{c}_1 \cdots \widehat{c}_r = \widehat{c} \neq 0$, where $l = (\sum^r_{i = 1} l_i - |J_i| - 1) + r - 1$, $l_i - |J_i| - 1 \geq 0$ and $J = J_1 \sqcup \cdots \sqcup J_r$. It follows that
\[
  n-1=\dim\K \geq l = \Bigl(\sum^r_{i = 1} l_i - |J_i| - 1\Bigr) + r - 1 \geq r - 1,
\]
hence $n \geq r$, as claimed.
\end{solution}

A (convex) \emph{polytope} $P$ is a bounded intersection of a finite number of halfspaces in a real affine space. A \textit{facet} of $P$ is its face of codimension~$1$.

A polytope $P$ of dimension $n$ is called \textit{simple} if each vertex of $P$ belongs to exactly $n$ facets. So if $P$ is simple, then the dual polytope $P^*$ is simplicial and its boundary $\6 P^*$ is a simplicial complex, which we denote by $\K_P$. Then $\K_P$ is the nerve complex of the covering of $\partial P$ by its facets. The moment-angle complex $\Z_{\K_P}$ is denoted simply by~$\Z_P$.

A \emph{simplicial sphere} (or \emph{triangulated sphere}) is a simplicial complex $\K$ whose geometric realisation is homeomorphic to a sphere. If $P$ is a simple polytope of dimension~$n$, then the nerve complex $\K_P$ is a simplicial sphere of dimension~$n-1$.
For $n\le3$, any simplicial sphere of dimension $n-1$ is combinatorially equivalent to the nerve complex $\K_P$ of a simple $n$-dimensional polytope~$P$. This is not true in dimensions $n\ge4$; the \emph{Barnette sphere} is a famous example of a $3$-dimensional simplicial sphere with $8$ vertices that is not combinatorially equivalent to the boundary of a convex $4$-dimensional polytope (see~\cite[\S2.5]{BP}).

\begin{thm}[{\cite[Theorem 4.1.4, Corollary 6.2.5]{BP}}]\label{dim2}
Let $\K$ be a simplicial sphere of dimension $(n - 1)$ with $m$ vertices. Then $\Z_\K$ is a closed topological manifold of dimension $m+n$. If $P$ be a simple $n$-dimensional polytope with $m$ facets, then $\Z_P$ is a smooth manifold of dimension $m+n$.
\end{thm}

A simple polytope $Q$ is called \textit{stacked} if it can be obtained from a simplex by a sequence of stellar subdivisions of facets. Equivalently, the dual simple polytope $P=Q^*$ is obtained from a simplex by iterating the vertex cut operation.

A \emph{connected sum of products of spheres} is a closed $n$-dimensional manifold $M$ homeomorphic to a connected sum $M_1\#\cdots\# M_k$ where each $M_k$ is a product of spheres $S^{n_{k1}}\times\cdots\times S^{n_{kl}}$, where $n_{k1}+\cdots+n_{kl}=n$.

The next theorem follows from the results of McGavran~\cite{mcga79}, see~\cite[Theorem~6.3]{bo-me06}. See also~\cite[\S2.2]{gi-lo13} for a different approach.

\begin{thm}[{see \cite[Theorem~4.6.12]{BP}}]\label{T4.6.12}
Let $P$ be a dual stacked $n$-polytope with $m > n + 1$ facets. Then the corresponding moment-angle manifold is homeomorphic to a connected sum of products of spheres with two spheres in each product, namely,
\[
\Z_P \cong \underset{k=3}{\overset{m-n+1}{\#}} (S^k \x S^{m+n-k})^{\#(k-2)\binom{m-n}{k-1}}
\]
\end{thm}
In particular, the moment-angle complex corresponding to a polygon (a two-dimensional polytope) is a connected sum of products of spheres.

A \emph{graph} $\Gamma$ is a one-dimensional simplicial complex.
A graph $\Gamma$ is called \emph{chordal} if every cycle of $\Gamma$ with more than $3$ vertices has a chord, where a chord is an edge connecting two vertices that are not adjacent in the cycle.
The vertices of a graph are in a \emph{perfect elimination order} if for any vertex $\{i\}$ all its neighbours with indices less than $i$ are pairwise adjacent.

\begin{thm}[\cite{FG}] \label{T_FG} 
A graph is chordal if and only if its vertices can be arranged in  a perfect elimination order. 
\end{thm}

The following property of chordal graphs is immediate from Theorem~\ref{T_FG}.


\begin{prop} \label{prop1}
Let $\Gamma$ be a chordal graph on $m$ vertices, and suppose that the vertices of $\Gamma$ are arranged in a perfect elimination order. Then $\Gamma \sm \{m\}$ is also a chordal graph, and the vertices of $\Gamma \sm \{m\}$ are automatically arranged in a perfect elimination order.
\end{prop}

\begin{lemma}\label{pairs}
Let $\K$ be a simplicial sphere of dimension greater than $1$ such that $H^*(\Z_{\K})\cong H^*(M_1 \# M_2 \# \cdots \# M_k)$ where each $M_i$ is a product of two spheres. Then the one-skeleton $\K^1$ is a chordal graph.
\end{lemma}
\begin{solution}
Let $\dim\K=n-1$ and $M_i = S^{l_i} \x S^{m+n-l_i}$, $i=1,\ldots,k$. We denote the corresponding generators of $H^*(\Z_\K)$ by $a_i$, $b_i$, where $\deg a_i=l_i$, $\deg b_i = m + n - l_i$, $i = 1, \ldots, k$, and $c$, $\deg c=m+n$ (the fundamental class). We have relations $a_i \cdot b_i = c$ for $i = 1, \ldots, k$, and all other products in $H^*(\Z_\K)$ are trivial. 

Suppose that there is a chordless cycle $C$ in $\K$ with $p > 3$ vertices. Then $C$ is a full subcomplex in $\K$, therefore $H^*(\Z_C)$ is a subring of $H^*(\Z_\K)$ by Lemma~\ref{retract}. By Theorem~\ref{T4.6.12} $\Z_C$ is also a connected sum of products of spheres, so there are nontrivial products $a'_j \cdot b'_j = c'$ in the ring $H^*(\Z_C)$, where $c'$ is the fundamental class of $\Z_C$ and $\deg c' = |C| + 2 \leq m + 2 < m + n = \deg c$, which is impossible in $H^*(\Z_\K)$. 
%
Thus, there are no chordless cycles in $\K$ with more than three vertices, so $\K^1$ is a chordal graph.
\end{solution}

The converse of Lemma~\ref{pairs} holds for two- and three-dimensional spheres, as shown in the next two sections, but fails in higher dimensions, as shown by the example below. A \emph{missing edge} of $\mathcal K$ is a pair of vertices that do not form a $1$-simplex.

\begin{expl}\label{5dimex}
Let $P$ be a three-dimensional polytope obtained by cutting two vertices of a tetrahedron~$\D^3$. By Theorem~\ref{T4.6.12},
\[
  \Z_P \cong (S^3 \x S^{6})^{\#3} \# (S^4 \x S^{5})^{\#2}  \,.
\]
Now let $P' = P \x \D^d$, where $d>1$, so that $\K_{P'}$ is a simplicial sphere of dimension $d+2>3$.
We have $\Z_{P'}=\Z_P\times \Z_{\D^d} \cong \Z_P \times S^{2d - 1}$, which is not a connected sum of products of spheres. However, $\K_{P'}^1$ is a chordal graph. Indeed, $\K_{P'} = \K_P * \6\D^d$. Hence, each missing edge of $\K_{P'}$ is a missing edge of~$\K_P$. There are only three missing edges in $\K_{P'}$, and no two of them form a chordless $4$-cycle. Also, there can be no chordless cycles with more than $4$ vertices, as any such chordless cycle has at least $5$ missing edges.
\end{expl}

The next lemma builds upon the results of~\cite[\S4]{FCMW}.

\begin{lemma}\label{corr_FCMW}
Let $\K$ be a simplicial sphere of dimension $>1$ such that 
$H^*(\Z_K) \cong H^*(M_1\# M_2\#\cdots\# M_k)$, where each $M_i$ is a product of spheres. Suppose that
$\K^1$ is not a chordal graph. Then all missing edges $I_1,\ldots,I_r$ of $\K$ are pairwise disjoint and
\[
  \K_{I_1\sqcup I_2\sqcup\cdots\sqcup I_r}=\K_{I_1} * \K_{I_2} * \cdots * \K_{I_r}.
\]
%
\end{lemma}
\begin{solution}
By \cite[Lemma~4.5]{FCMW} any chordless cycle in $\K^1$ has three or four vertices. Since $\K^1$ is not chordal, it contains a chordless $4$-cycle. Then by~\cite[Lemma~4.6]{FCMW} missing edges of $\K$ are pairwise disjoint, i.\,e. each pair of missing edges forms a chordless $4$-cycle. 

We have $H^3(\Z_\K) \cong \bigoplus_{|J| = 2}\wt{H}^0(\K_J) = \bigoplus_{j=1}^r \wt{H}^0(\K_{I_j})$ by Theorem~\ref{Hochster}. Choose a basis $a_1,\ldots,a_r$ of $H^3(\Z_\K)$ according to this decomposition, so that $a_j$ corresponds to a generator of $\widetilde H^0(\K_{I_j})=\widetilde H^0(S^0)\cong\mathbb Z$ for $j=1,\ldots,r$. Each product $a_j\cdot a_k$ is nonzero by Theorem~\ref{Hochster}, because $\K_{I_j\sqcup I_k}$ is a $4$-cycle.

Through the ring isomorphism $H^*(\Z_K) \cong H^*(M_1\# M_2\#\cdots\# M_k)$, three-\-di\-men\-si\-onal sphere factors $S^3_{ji}$ in the connected summands $M_i$ correspond to cohomology classes in $H^3(\Z_K)$, which we denote by $s_1,\ldots,s_r$. We have $H^3(\Z_K)\cong\mathbb Z\langle a_1,\ldots,a_r\rangle\cong\mathbb Z\langle s_1,\ldots,s_r\rangle$. Furthermore, if we denote the subring of $H^*(\Z_K)$ generated by $a_1,\ldots,a_r$ by $A$ and denote the subring generated by $s_1,\ldots,s_r$ by $R$, then we have a ring isomorphism $A\cong R$. Since $a_i\cdot a_j\ne 0$ for any $i\ne j$, we have $\mathop\mathrm{rank} A^6=\mathop\mathrm{rank} R^6=\frac{r(r-1)}2$. This implies that $s_i\cdot s_j\ne 0$ for $i\ne j$. It follows that all spheres $S^3_{ji}$, $j=1,\ldots,r$, belong to the same connected summand~$M_i$, because the product of the cohomology classes corresponding to sphere factors in different summands of the connected sum $M_1\# M_2\#\cdots\# M_k$ is zero. Therefore, $s_1\cdot s_2\cdots s_r\ne0$ in~$R$. This implies, by the ring isomorphism $A\cong R$, that $a_1\cdot a_2\cdots a_r$ is nonzero in $H^*(\Z_\K)$. Now it follows from the product description in Theorem~\ref{Hochster} that $\K_{I_1\sqcup I_2\sqcup\cdots\sqcup I_r}=\K_{I_1} * \K_{I_2} * \cdots * \K_{I_r}$.
\end{solution}

\section{Two-dimensional spheres}
Here we consider moment-angle manifolds corresponding to two-dimensional simplicial spheres $\K$ or, equivalently, to three-dimensional simple polytopes $P$. 

The case $P=I^3$ (a three-dimensional cube) is special. In this case the nerve complex $\K_P$ is $S^0*S^0*S^0$ (the join of three $0$-dimensional spheres, or the boundary of a three-dimensional cross-polytope) and $\Z_P\cong S^3\times S^3\times S^3$.

A simplicial complex $\K$ is called \textit{Golod} if the multiplication and all higher Massey products in $H^*(\Z_\K)$ are trivial. (Equivalently, the Stanley--Reisner ring $\k[\K]$ is a Golod ring over any field $\k$, see~\cite[\S4.9]{BP}.) A simplicial complex $\K$ on $[m]$ is called \textit{minimally non-Golod} if $\K$ is not Golod, but for any vertex $i\in [m]$ the complex $\K_{[m]\setminus\{i\}}$ is Golod.

The following result extends~\cite[Proposition~11.6]{bo-me06}, where the equvalence of conditions (a), (b) and (c) was proved:

\begin{prop}\label{mainthm3}
Let $\K$ be a two-dimensional simplicial sphere and let $P$ be the a three-dimensional simple polytope such that $\K=\K_P$. Suppose that $P$ is not a cube. The following conditions are equivalent:
\begin{itemize}
    \item[(a)] $P$ is obtained from a simplex $\D^3$ by iterating the vertex cut operation, i.\,e. $P^*$ is a stacked polytope;
    \item[(b)] $\Z_P$ is diffeomorphic to a connected sum of products of spheres;
    \item[(c)] $H^*(\Z_P)$ is isomorphic to the cohomology ring of a connected sum of products of spheres;
    \item[(d)] the one-dimensional skeleton of the nerve complex $\K_P$ is a chordal graph;
    \item[(e)] $\K_P$ is minimally non-Golod, unless $P=\Delta^3$.
\end{itemize}
\end{prop}


\begin{solution}
We prove the implications (a)$\Rightarrow$(b)$\Rightarrow$(c)$\Rightarrow$(d)$\Rightarrow$(a), (e)$\Rightarrow$(d) and (a)$\Rightarrow$(e).

\smallskip

\noindent (a)$\Rightarrow$(b) This is Theorem~\ref{T4.6.12}.

\smallskip

\noindent (b)$\Rightarrow$(c) is clear.

\smallskip

\noindent (c)$\Rightarrow$(d) Let $H^*(\Z_P) \cong H^*(M_1 \# M_2 \# \cdots \# M_k)$, where each $M_i$ is a product of spheres. 
Since the cohomological product length of $\Z_P$ is at most $3$ 
(Proposition~\ref{cohom_l}), there is at most $3$ sphere factors in each~$M_i$. If some $M_i$ has exactly $3$ factors, then $\Z_P=S^3\times S^3\times S^3$ and $P$ is a cube by~\cite[Theorem~4.3~(a)]{FCMW}. This contradicts the assumption. 
Now, $\K^1_P$ is a chordal graph by Lemma~\ref{pairs}.

\smallskip

\noindent (d)$\Rightarrow$(a) We use induction on the number $m$ of facets of $P$. The base $m = 4$ is clear, as $P$ is a simplex $\D^3$ in this case.

For the induction step, assume that the vertices of $\K_P$ are arranged in a perfect elimination order. Let $j_1, \ldots, j_s$ be the vertices adjacent to the last vertex~$m$. First we prove that $s = 3$. 

Let $F_i$ denote the $i$th facet of~$P$. Since $\{j_1, \ldots, j_s\}$ is a clique of $\K^1_P$, the facets $F_{j_1},\ldots,F_{j_s}$ are pairwise adjacent.
Suppose that $s \geq 4$. Renumbering the facets if necessary, we may assume that $F_{j_1}$, $F_{j_2}$, $F_{j_3}$, $F_{j_4}$ are consecutive facets in a cyclic order around~$F_m$, so that $F_m \cap F_{j_1} \cap F_{j_3} = \varnothing$ and $F_m \cap F_{j_2} \cap F_{j_4} = \varnothing$. 
Since $F_{j_1}$ and $F_{j_3}$ are adjacent, the facets $F_m$, $F_{j_1}$ and $F_{j_3}$ form a $3$-belt (a prismatic $3$-circuit). This $3$-belt splits $\partial P$ into two connected components~\cite[Lemma~2.7.2]{E}. The facets $F_{j_2}$ and $F_{j_4}$ lie in different components, so they cannot be adjacent. A contradiction. Hence, $s=3$.

Since $F_m$ has $3$ adjacent facets, it is a triangle. If $F_m$ is adjacent to a triangular facet, then $P$ is a simplex. Otherwise, there exist a polytope $P'$ such that $P$ is obtained from $P'$ by cutting a vertex with formation of a new facet~$F_m$. 
Then $\K_{P'}$ is obtained from $\K_P$ by removing the vertex $\{m\}$ and adding simplex $\{j_1, j_2, j_3\}$. Hence, the $1$-skeleton of $\K_{P'}$ is also a chordal graph by Proposition~\ref{prop1}. Now $P'$ has $m-1$ facets, so we complete the induction step.

\smallskip

(e)$\Rightarrow$(d) Let $\K_P$ be minimally non-Golod, and suppose there is a chordless cycle $C$ in $\K_P^1$ with $p > 3$ vertices. Then $C$ is a full subcomplex of $\K_P$ and $p < m$ (otherwise $\K_P = C$, which is impossible for a $3$-dimensional polytope). For any vertex $v \in [m] \sm C$, note that $C$ is also a full subcomplex in $\K_P \sm \{v\}$. Therefore, $H^*(\Z_C)$ is a subring of $H^*(\Z_{\K_P \sm \{v\}})$ by Lemma~\ref{retract}. On the other hand, there are nontrivial products is $H^*(\Z_C)$ by Theorem~\ref{T4.6.12}, whereas all products in $H^*(\Z_{\K_P \sm \{v\}})$ must be trivial, since $\K_P \sm \{v\}$ is Golod. A contradiction. Hence, there are no chordless cycles in~$\K_P^1$.  

\smallskip

(a)$\Rightarrow$(e) This follows from~\cite[Theorem~3.9]{L}: if an $n$-dimensional simple polytope $P$ is obtained from $P'$ by a vertex cut, and $\K_{P'}$ is minimally non-Golod, then $\K_P$ is also minimally non-Golod.
\end{solution}

\section{Three-dimensional spheres}
Recall that the product in $H^*(\Z_\K)=\H^{*,*}(\K)$ is given by~\eqref{multH}.
A nonzero element $c \in \mathcal H^{l,J}= \wt{H}^l(\K_J)$ is \textit{decomposable} if $c = \sum_{i=1}^p a_i \cdot b_i$ for some nonzero $a_i \in \wt{H}^{r_i}(\K_{I_i})$, $b_i \in \wt{H}^{l-1-r_i}(\K_{J \sm I_i})$, where $0 \leq r_i \leq l - 1$ and $I_i \ss J$ are proper subsets for $i=1,\ldots,p$.


A \emph{missing face} (or a \emph{minimal non-face}) of $\K$ is a subset $I\subset[m]$ such that $I$ is not a simplex of~$\K$, but every proper subset of $I$ is a simplex of~$\K$. Each missing face corresponds to a full subcomplex $\partial\Delta_I\subset\K$, where $\partial\Delta_I$ denotes the boundary of simplex $\Delta_I$ on the vertex set~$I$. A missing face $I$ defines a simplicial homology class in $\widetilde H_{|I|-2}(\K)$, which we continue to denote by~$\partial\Delta_I$.  
We denote by $\mathop\mathrm{MF}_n(\K)$ the set of missing faces $I$ of dimension $n$, that is, with $|I|=n+1$.


\begin{lemma}\label{MFs}
Let $I\in \mathop\mathrm{MF}_l(\K)$ be a missing face of $\K$. Then any cohomology class $c \in \H^{l-1, *}(\K)$ such that $\langle c,\6\D_I\rangle \neq 0$ is indecomposable.
\end{lemma}
\begin{solution}
Let $\K'$ be the simplicial complex obtained from $\K$ by filling in all missing faces of dimension~$l$ with simplices, so that $\mathop\mathrm{MF}_l(\K') = \varnothing$ and $\K^{l-1}=(\K')^{l-1}$. Then the inclusion $i: \K \hookrightarrow \K'$ induces a ring homomorphism $i^*: \H^{*,*}(\K') \rightarrow \H^{*,*}(\K)$ and $\H^{r,*}(\K') \cong \H^{r,*}(\K)$ for $r \leq l-2$. Also, $i_*(\6\D_I) = 0$ for $I\in\mathop\mathrm{MF}_l(\K)$.

Suppose $c$ is decomposable, that is, $c = \sum_{i=1}^p a_i \cdot b_i$. Choose $a_i', b_i'$ such that $i^*(a_i') = a_i$ and $i^*(b_i') = b_i$ and define $c' := \sum_{i=1}^p a_i' \cdot b_i'$. Then $i^*(c') = c$ and
\[
  \langle c,\6\D_I\rangle = \langle i^*(c'),\6\D_I\rangle = 
  \langle c',i_*(\6\D_I)\rangle = 0.
\]
This is a contradiction.
\end{solution}






\begin{thm}\label{chordal4}
Let $\K$ be a three-dimensional simplicial sphere such that $\K\ne\partial\Delta^4$ and $\K^1$ is a chordal graph. Then $H^*(\Z_\K) \cong H^*(M)$, where $M$ is a connected sum of products of spheres with two spheres in each product.
\end{thm}
\begin{solution}
We use the notation $\H^{*,*}=H^*(\Z_\K)$ and analyse possible nontrivial products in~\eqref{multH}. We have $\H^{k,*}=0$ for $k\geq4$ since $\K$ is a three-dimensional sphere.
Products of the form $\H^{3,*} \ox \H^{i,*} \to \H^{4+i, *}$, $\H^{2,*} \ox \H^{2,*} \to \H^{5,*}$ and $\H^{2,*} \ox \H^{1,*} \to \H^{4,*}$ are therefore trivial for dimensional reasons. 

Since $\K$ is a 3-dimensional sphere, $\Z_\K$ is an $(m+4)$-dimensional manifold. Nontrivial products $\wt{H}^i(\K_I) \ox \wt{H}^{2-i}(\K_J) \to \wt{H}^3(\K_{I \cup J})$ come from Poincar\'e duality for~$\Z_\K$ (see~\cite[Proposition~4.6.6]{BP}), because $\wt{H}^3(\K_{I \cup J})$ is nonzero only when $I\sqcup J=[m]$. The Poincar\'e duality isomorphisms $\wt{H}^i(\K_I) \cong \wt{H}_{2-i}(\K_{[m] \sm I})$ (or the Alexander duality isomorphisms for the $3$-sphere~$\K$, see~\cite[3.4.11]{BP}) imply that the groups $\wt{H}^i(\K_I)$ are torsion-free for any $i$ and $I \ss [m]$.

Next we prove that all multiplications of the form $\H^{0,*} \ox \H^{0,*} \longrightarrow \H^{1,*}$ are trivial. 
Assume that there are cohomology classes $a, b \in \H^{0, *}$ such that $0 \neq a \cdot b =: c \in \wt{H}^1(\K_I)$. Since $c \neq 0$ there exists $\gamma \in H_1(\K_I)$ such that $\langle c, \gamma\rangle \neq 0$. We can write $\gamma = \lambda_1 \gamma_1 + \cdots + \lambda_k\gamma_k$, where each $\gamma_i$ is a simple chordless cycle in $\K^1$ and $\lambda_i\ne0$. Since  $\K^1$ is chordal, $\gamma_i \in \mathop\mathrm{MF}_2(\K)$. Now, $0 \neq \langle c, \gamma \rangle = \sum_{j=1}^k \lambda_i \langle c, \gamma_i\rangle$, so $\langle c,\gamma_i\rangle \neq 0$ for some~$i$. Hence, $c$ is indecomposable by Lemma~\ref{MFs}. A contradiction.

Finally, we prove that all multiplications of the form $\H^{0,*} \ox \H^{1,*} \longrightarrow \H^{2,*}$ are trivial.
Assume that there exists a nontrivial product $a^0 \cdot b^1 = c^2 \neq 0$ for some $a^0 \in \wt{H}^0(\K_I)$, $b^1 \in \wt{H}^1(K_J)$, $c^2 \in \wt{H}^2(\K_{I \cup J})$. By Poincar\'e duality there exists an element $a' \in \wt{H}^0(\K_{[m] \sm (I \cup J)})$ such that $0 \neq a' \cdot c^2 = a' \cdot a^0 \cdot b^1 \in \wt{H}^3(\K)$. Then $a^0 \cdot a' \neq 0$, so we obtain a nontrivial multiplication of the form $\H^{0,*} \ox \H^{0,*} \longrightarrow \H^{1,*}$. A contradiction.

It follows that the only nontrivial multiplications in $\H^{*,*}(\K)$ are 
\[
  \H^{0,I} \ox \H^{2,[m]\setminus I} \longrightarrow \H^{3,[m]}\quad 
  \text{and}\quad 
  \H^{1,J} \ox \H^{1,[m]\setminus J} \longrightarrow \H^{3,[m]},
\]  
which arise from Poincar\'e duality. Therefore, the ring $H^*(\Z_\K)$ is free as an abelian group with $\mathbb Z$-basis 
\[
  \{1,a^0_1,\ldots,a^0_k,a^1_1,\ldots,a^1_l,b^1_1,\ldots,b^1_l,
  b^2_1,\ldots,b^2_k,c\},
\]  
where $a^0_1,\ldots,a^0_k\in\H^{0,*}$, $a^1_1,\ldots,a^1_l,b^1_1,\ldots,b^1_l\in\H^{1,*}$,
$b^2_1,\ldots,b^2_k\in\H^{2,*}$, $c\in\H^{3,m}=H^{m+3}(\Z_\K)$ is the fundamental class, and the product is given by
$a^0_i \cdot b^2_j = \delta_{ij}c$ and $a^1_p \cdot b^1_q = \delta_{pq}c$, where $\delta_{ij}$ is the Kronecker delta. At least one of the groups $\H^{0,*}$ and $\H^{1,*}$ is nonzero, as otherwise $\mathcal K=\partial\Delta^4$ and $\Z_\K\cong S^9$. Then $H^*(\Z_\K)=\H^{*,*}$ is isomorphic to the cohomology ring of a connected sum of products of spheres with two spheres in each product.
\end{solution}

For simplicial spheres $\K$ of dimension $>3$, the condition that $\mathcal K^1$ is a chordal graph does not imply that $H^*(\mathcal Z_{\mathcal K})$ is isomorphic to the cohomology ring of a connected sum of products of spheres, as shown by Example~\ref{5dimex}. The next result gives a sufficient condition in any dimension. We say that the group $\H^{l,*}(\K)$ is \emph{generated by missing faces} of $\K$ if for any nonzero $c \in \H^{l,*}(\K)$ there exists $I \in \mathop\mathrm{MF}_{l+1}(\K)$ such that $\langle c, \6\D_I \rangle \neq 0$.

\begin{thm}\label{ssndim}
Let $\K$ be a simplicial sphere of dimension $d$ such that $\K\ne\partial\Delta^{d+1}$ and the group $\H^{l,*}(\K)$ is generated by missing faces of $\K$ for $l \leq \left\lfloor\frac{2d-1}{3}\right\rfloor$. Then $H^*(\Z_\K)$ is isomorphic to the cohomology ring of a connected sum of products of spheres with two spheres in each product.
\end{thm}
\begin{solution}
We can assume that $d\ge2$, as otherwise $\mathcal K$ is the boundary of polygon and the result follows from Theorem~\ref{T4.6.12}.
As in the proof of Theorem~\ref{chordal4}, we analyse possible nontrivial products in~\eqref{multH}. 
We denote $q := \left\lfloor\frac{2d-1}{3}\right\rfloor$.

We have $\H^{k,*}=0$ for $k>d$ since $\K$ is an $d$-dimensional sphere.
Therefore, products of the form $\H^{i,*} \ox \H^{j,*} \to \H^{i+j+1, *}$ with $i+j \ge d$ are trivial. 

Nontrivial products or the form $\H^{i,*} \ox \H^{j,*} \to \H^{i+j+1, *}$ with $i+j = d-1$ are given by $\wt{H}^i(\K_I) \ox \wt{H}^{d-1-i}(\K_J) \to \wt{H}^d(\K_{I \cup J})$ and come from Poincar\'e duality, because $\wt{H}^d(\K_{I \cup J})$ is nonzero only when $I\sqcup J=[m]$. 
We prove by contradiction that the groups $\wt{H}^i(\K_I)$ are torsion-free for $i \leq q$. Assume that there is a cocycle $0 \neq c \in \H^{i,*}(\K)$ and a nonzero integer $k$ such that $k \cdot c = 0$. Let $\Tilde{c}$ be a representing cochain for $c$, then $k \cdot \Tilde{c}$ is a coboundary and $k\cdot \Tilde{c} = d\Tilde{b}$ for some cochain $\Tilde{b}$. By assumption there exists $I \in \mathop\mathrm{MF}_{i+1}(\K)$ such that $\langle c, \6\D_I \rangle \neq 0$, hence,
\[
0 \neq k \cdot \langle c, \6\D_I \rangle = \langle k \cdot \Tilde{c}, \6\D_I \rangle = \langle d\Tilde{b}, \6\D_I \rangle = \langle \Tilde{b}, \6(\6\D_I) \rangle = 0
\]
and we get a contradiction. Now the Alexander duality isomorphisms $\widetilde H^i(\K_J) \cong \wt{H}_{d-1-i}(\K_{[m] \sm J})$ imply that the homology groups $\wt{H}_j(\K_J)$ are torsion-free for $j\ge d-1-q$. Since $d-1-q\le q$, we obtain that $\wt{H}_j(\K_J)$ is torsion-free for $j\ge q$, whereas $\wt{H}^j(\K_J)$ is torsion-free for $j\le q$. By the universal coefficient theorem we conclude that the groups $\wt{H}^j(\K_J)$ are torsion-free for all $j$ and~$J$.

All products of the form $\H^{i,*} \ox \H^{j,*} \longrightarrow \H^{i+j+1,*}$ are trivial for $i+j < q$, since any $l$-dimensional cohomology class with $l \leq q$ is indecomposable by Lemma~\ref{MFs}.

Finally, we prove that all products of the form $\H^{i,*} \ox \H^{j,*} \longrightarrow \H^{i+j+1,*}$ are trivial for $q \leq i+j \leq d - 2$. Suppose there are classes $a \in \H^{i, I}$ and $b \in \H^{j, J}$ with $q \leq i+j \leq d - 2$ such that $0 \neq a \cdot b =: c \in \H^{i+j+1, I\cup J}$. Without loss of generality we assume that $i \leq j$. Then there exists an element $a' \in \wt{H}^{d-i-j-2}(\K_{[m] \sm (I \cup J)})$ such that $0 \neq a' \cdot c = a' \cdot a \cdot b \in \wt{H}^d(\K)$ by Poincar\'e duality. Therefore, $a \cdot a' \neq 0$ and so we obtain a nontrivial product of the form $\H^{i,*} \ox \H^{k,*} \longrightarrow \H^{i+k+1,*}$ for $k = d-i-j-2$. By assumption, $q \leq i + j \leq 2j$ and $q>\frac{2d-1}3-1$, hence,
\[
  i + k = d-j-2 \le d-2-\frac{q}{2} < q.
\]
Thus, $a' \cdot a$ is a product of the form $\H^{i,*} \ox \H^{k,*} \longrightarrow \H^{i+k+1,*}$ with $i+k < q$, so it must be trivial. A contradiction.

We obtain that the only nontrivial products in $\H^{*,*}(\K)$ arise from Poincar\'e duality. It follows that the ring $H^*(\Z_\K)$ is isomorphic to the cohomology ring of a connected sum of products of spheres with two spheres in each product.
\end{solution}

The next theorem extends the result of Theorem~\ref{chordal4} to a complete characterisation of three-dimensional spheres $\K$ such that $H^*(\Z_\K)$ is isomorphic to the cohomology ring of a connected sum of products of spheres.

\begin{thm}\label{chordal4iff}
Let $\K$ be a three-dimensional simplicial sphere. Then $H^*(\Z_\K) \cong H^*(M_1\#\cdots\# M_k)$ where each $M_i$ is a product of spheres if and only if one of the following conditions is satisfied:
\begin{itemize}
\item[(a)] $\K=S^0*S^0*S^0*S^0$ (the boundary of a $4$-dimensional cross-polytope);
\item[(b)] $\K^1$ is a chordal graph;
\item[(c)] $\K^1$ has exactly two missing edges which form a chordless $4$-cycle.
\end{itemize}
\end{thm}

\begin{solution}
First we prove the ``only if'' statement. If $\K^1$ is a chordal graph, then (b) is satisfied. Otherwise, by Lemma~\ref{corr_FCMW} the missing edges $I_1, \ldots, I_r$ of $\K$ are pairwise disjoint and $\K_{I_1 \sqcup \cdots \sqcup I_r}=\K_{I_1} * \cdots *\K_{I_r}$. We have $r\le 4$, since $\dim\K=3$. 

If $r=4$, then $\K=\K_{I_1} * \cdots *\K_{I_4}$, so that (a) holds. 

If $r = 3$, then $\K_{I_1 \sqcup I_2 \sqcup I_3}=\K_{I_1} * \K_{I_2} *\K_{I_3}$ is a two-dimensional simplicial sphere. We have $\wt{H}_0(\K \sm \K_{I_1 \sqcup I_2 \sqcup I_3}) \cong \wt{H}^2(\K_{I_1 \sqcup I_2 \sqcup I_3}) \cong \mb{Z}$ by Alexander duality. Hence, $\K \sm \K_{I_1 \sqcup I_2 \sqcup I_3}$ is not connected. It follows that there is at least one more missing edge in $\K$ besides $I_1, I_2, I_3$. A contradiction.

If $r = 2$, then (c) holds.
    
If $r = 1$, then $\K^1$ is in fact a chordal graph, since any chordless cycle with more than three vertices has at least two missing edges. Hence, (b) holds.

\smallskip

Now we prove the ``if'' statement. If (a) holds, then $\mathcal Z_\K$ is a product of spheres. If (b) holds, then $H^*(\Z_\K) \cong H^*(M_1\#\cdots\# M_k)$ where each $M_i$ is a product of spheres by Theorem~\ref{chordal4}. Suppose (c) holds. Then $\H^{0,*}(\K) = \mb{Z}\langle a_1, a_2 \rangle$, where $a_1$ and $a_2$ correspond to the two missing edges of $\K$, and $a_1 \cdot a_2 \neq 0$. We use the same argument as in the proof of Theorem~\ref{chordal4} with one exception: there is one nontrivial product of the form $\H^{0,*}(\K) \ox \H^{0,*}(\K) \ox \H^{1,*}(\K) \longrightarrow \H^{3,*}(\K)$. Namely, $a_1 \cdot a_2 \cdot b \mapsto c$, where $b$ is Poincar\'e dual to $a_1 \cdot a_2$ and $c$ is the fundamental class of~$\K$. All other nontrivial products in $H^*(\Z_\K)$ arise from Poincar\'e duality. Thus the ring $H^*(\Z_\K)$ is generated by elements $\{a_1, a_2, b, c, x_i, y_i \colon i = 1, 2, \ldots, N\}$, where $x_i, y_i\in\H^{1,*}(\K)$, with the following multiplication rules: $a_1 \cdot a_2 \cdot b = c$, $x_i \cdot y_i = c$ for $i = 1, 2, \ldots, N$, and all other products of generators are zero. Clearly, $H^*(\Z_\K)$ is isomorphic to the cohomology ring of a connected sum of products of spheres.
\end{solution}

\begin{rem}
Note that under condition (c) of Theorem~\ref{chordal4iff} we have $H^*(\Z_\K) \cong H^*(M)$, where $M$ is a connected sum of products of spheres in which one of the summands is a product of \emph{three} spheres. The first example of such a simplicial sphere $\K$  was constructed in~\cite{FCMW}. Later it was shown in~\cite{iriy18} that the corresponding moment-angle manifold $\mathcal Z_\K$ is diffeomorphic to~$M$. 
\end{rem}

\begin{rem}
It can be shown that if $\K$ is a three-dimensional simplicial sphere such that $\K^1$ is a chordal graph, then all higher Massey products in $H^*(\Z_\K)$ are trivial. This implies that a three-dimensional simplicial sphere $\K\ne\partial\Delta^4$ is minimally non-Golod if and only if $\K^1$ is a chordal graph. We elaborate on this in a subsequent paper.
\end{rem}

%
%
%

\end{document}